\documentclass[12pt]{amsart}
\usepackage{amsmath,amscd,amssymb,amsfonts,graphics}
\usepackage[bbgreekl]{mathbbol}
\DeclareSymbolFontAlphabet{\mathbbold}{bbold}
\newcommand{\h}{\hbox}
\newcommand{\q}{\quad}
\newcommand{\nin}{\noindent}
\newcommand{\bs}{\par\bigskip}
\newcommand{\ms}{\par\medskip}
\newcommand{\sk}{\par\smallskip}
\newcommand{\bsn}{\par\bigskip\noindent}
\newcommand{\msn}{\par\medskip\noindent}
\newcommand{\skn}{\par\smallskip\noindent}
\newcommand{\ges}{\geqslant}
\newcommand{\les}{\leqslant}
\newcommand{\1}{\hskip1pt}
\newcommand{\mopl}{\hbox{$\bigoplus$}}
\newcommand{\mcap}{\hbox{$\bigcap$}}

\newcommand{\msum}{\hbox{$\sum$}}

\newcommand{\D}{{\mathcal D}}

\newcommand{\Hc}{{\mathcal H}}
\newcommand{\Lc}{{\mathcal L}}
\newcommand{\M}{{\mathcal M}}

\newcommand{\DD}{{\mathbf D}}

\newcommand{\R}{{\mathbb R}}
\newcommand{\RR}{{\mathbf R}}
\newcommand{\Q}{{\mathbb Q}}
\newcommand{\C}{{\mathbb C}}
\newcommand{\Z}{{\mathbb Z}}

\newcommand{\dd}{\partial}
\newcommand{\ddd}{{\rm d}}
\newcommand{\Gr}{{\rm Gr}}
\newcommand{\Om}{\Omega}

\newcommand{\ssb}{\raise.15ex\h{${\scriptscriptstyle\bullet}$}}
\newcommand{\ssc}{\,\raise.15ex\hbox{${\scriptstyle\circ}$}\,}
\newcommand{\onto}{\mathop{\rlap{$\to$}\hskip2pt\hbox{$\to$}}}
\newcommand{\into}{\hookrightarrow}
\newcommand{\simto}{\,\,\rlap{\hskip1.3mm\raise1.4mm\hbox{$\sim$}}\hbox{$\longrightarrow$}\,\,}
\newcommand{\pl}{\1{+}\1 }
\newcommand{\mi}{\1{-}\1}
\newcommand{\eq}{\,{=}\,}
\newcommand{\nes}{\,{\ne}\,}
\newcommand{\sgt}{\,{>}\,}

\newcommand{\sst}{\,{\subset}\,}

\newcommand{\ins}{\,{\in}\,}
\newcommand{\tos}{\,{\to}\,}
\newcommand{\defs}{\,{:=}\,}

\begin{document}
\h{}\bs
\centerline{\large Local and Global Invariant Cycle Theorems for Hodge Modules}
\bs
\centerline{Morihiko Saito}
\ms\bsn
\vbox{\narrower\smaller\nin{\bf Abstract.} We show that the local and global invariant cycle theorems for Hodge modules follow easily from the general theory. We also give some comments about related papers.}
\ms\bs
\centerline{\bf Introduction}
\bsn
It does not seem well recognized (see for instance \cite{ES}) that the local and global invariant cycle theorems for pure Hodge modules follow easily from the general theory \cite{mhp}, \cite{mhm}. In these notes, we show that the decomposition theorem implies the {\it local invariant cycle theorem\1} for pure Hodge modules (see {\bf 1.1} below), and the {\it global invariant cycle theorem\1} for pure Hodge modules can be proved in a similar way to the classical case \cite[4.1.1\,(ii)]{De1}, see {\bf 1.2} below.
\sk
As for the estimate of weights of the cohomology of a link (which is called``local purity" in \cite{ES}), this has been known in the constant coefficient case (see \cite[1.18]{intr}, \cite{DS}), and a similar reasoning apply to the pure Hodge module case, since the assertion was proved using mixed Hodge modules, see {\bf 2.1} below.
\sk
This work was partially supported by JSPS Kakenhi 15K04816.
\bs\bs
\centerline{\bf 1. Local and global invariant cycle theorems.}
\bsn
{\bf 1.1.~Local invariant cycle theorem.} Let $f\,{:}\,X\to\Delta$ be a proper morphism from a complex manifold to a disk. Here we assume either $f$ is projective or $X$ is an open subset of a smooth complex algebraic variety. Let $\M$ be a pure Hodge module with strict support $Y$ which is not contained in a fiber of $f$. Let $K$ be the underlying $\Q$-complex of $\M$. Then in the notation of \cite{BBD}, we have the decomposition theorem asserting the non-canonical isomorphism
$$\RR f_*K\cong\mopl_k\,{}^{\bf p}\!R^k\!f_*K\q\q\h{with}\q\q{}^{\bf p}\!R^kf_*={}^{\bf p}\Hc^k\RR f_*,
\leqno(1.1.1)$$
together with the isomorphisms
$${}^{\bf p}\!R^k\!f_*K=(j_*L_{\Delta^*}^k)[1]\oplus L_0^k\q\q(k\in\Z).
\leqno(1.1.2)$$
Here $L^k_{\Delta^*}$, $L_0^k$ are local systems on $\Delta^*$, $0$, and $j:\Delta^*\into\Delta$ denotes the canonical inclusion. (This assertion can be reduced to the $f$ projective case.)
\sk
These isomorphisms give the non-canonical isomorphisms
$$R^k\!f_*K\cong j_*L_{\Delta^*}^{k+1}\oplus L_0^k\q\q(k\in\Z).
\leqno(1.1.3)$$
These imply the following.
\msn
{\bf Theorem~1.1} ({\it Local invariant cycle theorem}\1). {\it We have canonical surjection
$$H^k(X_0,K|_{X_0})\onto H^k(X_s,K|_{X_s})^T\q\q(s\in\Delta^*,\,k\in\Z),
\leqno(1.1.4)$$
shrinking $\Delta$ if necessary, where the right-hand side denotes the $T$-invariant subspace with $T$ the local monodromy.}
\msn
{\it Proof.} By the proper base change theorem, we have the isomorphisms
$$H^k(X_s,K|_{X_s})=(R^k\!f_*K)_s\q\q(s\in\Delta,\,k\in\Z).
\leqno(1.1.5)$$
So the assertion follows from (1.1.3). (Note that (1.1.4) is a property of the {\it sheaf\1} $R^k\!f_*K$, which depends only on the isomorphism class of the sheaf.)
\msn
{\bf 1.2.~Global invariant cycle theorem.} One can generalize an argument in \cite[4.1.1\,(ii)]{De1} as follows. Let $f:X\to S$ be a proper surjective morphism of irreducible complex algebraic varieties. Let $\M$ be a pure Hodge module of weight $w$ with strict support $X$, and $K$ be the underlying $\Q$-complex. We have the following.
\msn
{\bf Theorem~1.2} ({\it Global invariant cycle theorem}\1). {\it There is the canonical surjection for $s\in S'\,{:}$
$$H^k(X,K)\onto H^k(X_s,K|_{X_s})^{G_{k,s}}\q\q(k\in\Z).
\leqno(1.2.1)$$
Here $S'\subset S$ is a sufficiently small non-empty smooth Zariski-open subset such that the $R^k\!f_*K|_{S'}$ are local systems $(k\in\Z)$, and the $G_{k,s}$ denote the {\it monodromy group\1} of the local system $R^k\!f_*K|_{S'}$ with base point $s$.}
\msn
{\it Proof.} Set $X':=f^{-1}(S')$. Let $f':X'\to S'$ be the restriction of $f$. The decomposition theorem for $f'$ implies the canonical surjection
$$\Gr^W_{w+k}H^k(X',K|_{X'})\onto H^k(X_s,K|_{X_s})^{G_{k,s}}\q\q(s\in S',\,k\in\Z),
\leqno(1.2.2)$$
since the $R^k\!f_*K|_{S'}$ are local systems. Here $H^k(X_s,K|_{X_s})$ is pure of weight $w\pl k$. Indeed, $\M[-d_S]|_{X_s}$ is a pure Hodge module of weight $w\mi d_S$ on $X_s$ ($s\in S'$), and
$$H^k(X_s,K|_{X_s})=H^{k+d_S}(X_s,K[-d_S]|_{X_s})\q\q(d_S\,{:=}\,\dim S).$$
\sk
We then get (1.2.1) from (1.2.2), since we have moreover the canonical surjection
$$H^k(X,K)\onto\Gr^W_{w+k}H^k(X',K|_{X'}).
\leqno(1.2.3)$$
This surjection follows from the long exact sequence of mixed Hodge structures
$$\to H^k(X,K)\to H^k(X',K|_{X'})\to H^{k+1}(X'',i^!K)\to
\leqno(1.2.4)$$
with $X'':=X\,{\setminus}\,X'$ and $i:X''\into X$ the natural inclusion. Indeed, $H^{k+1}(X'',i^!K)$ has weights $\ges w\pl k\pl 1$, since $i^!\M$ has weights $\ges w$, see \cite[(4.5.2)]{mhm}. So Thm.\,1.2 follows.
\bs\bs
\centerline{\bf 2. Local purity in the sense of \cite{ES}.}
\bsn
{\bf 2.1.~Local purity.} Let $\M$ be a pure Hodge module of weight $w$ with strict support $X$. Take $x\in X$ with inclusions $i_x:\{x\}\into X$, $j_x:X\,{\setminus}\,\{x\}\into X$. Then the ``local purity" in the sense \cite{ES} asserts the following.
\msn
{\bf Theorem~2.1.} 
$$\h{\it $H^ki_x^*(j_x)_*j_x^*\M$ has weights $\les w\pl k\,$ if $\,k<0$, and $>w\pl k\,$ if $\,k\ges 0$.}
\leqno(2.1.1)$$
\skn
{\bf Remark~2.1a.} This is known in the constant coefficient case, see \cite[1.18]{intr}, \cite{DS}, where mixed Hodge modules are used for the proof. It is easy to generalize this as follows.
\msn
{\bf Proof of Theorem~2.1.}
Applying $i_x^*$ to the distinguished triangle
$$(i_x)_*i_x^!\M\to\M\to (j_x)_*j_x^*\M\buildrel{+1}\over\to,$$
\vskip-3mm\nin
we get
$$i_x^!\M\to i_x^*\M\to i_x^*(j_x)_*j_x^*\M\buildrel{+1}\over\to.
\leqno(2.1.2)$$
Taking its dual, and using the self-duality $\DD\M=\M(w)$, it gives
$$\DD\1i_x^*(j_x)_*j_x^*\M\to i_x^!\M(w)\to i_x^*\M(w)\buildrel{+1}\over\to,
\leqno(2.1.3)$$
since $\DD\1i_x^*=i_x^!\1\DD$. We thus get the self-duality
$$\DD\1i_x^*(j_x)_*j_x^*\M=i_x^*(j_x)_*j_x^*\M(w)[-1].
\leqno(2.1.4)$$
Setting $H^k:=H^ki_x^*(j_x)_*j_x^*\M$, this means the duality of mixed Hodge structures
$$\DD H^k=H^{-k-1}(w)\q\q(k\in\Z).
\leqno(2.1.5)$$
So the assertion (2.1.1) is reduced to the case $k<0$.
\sk
Consider the composition
$$(j_x)_*j_x^*\M\to(i_x)_*i_x^*(j_x)_*j_x^*\M\to\tau^{\ges 0}(i_x)_*i_x^*(j_x)_*j_x^*\M,
\leqno(2.1.6)$$
Let $\M''$ be its shifted mapping cone so that we have the distinguished triangle
$$\M''\to(j_x)_*j_x^*\M\to\tau^{\ges 0}(i_x)_*i_x^*(j_x)_*j_x^*\M\buildrel{+1}\over\to,
\leqno(2.1.7)$$
Let $K''$ be the underlying $\Q$-complex of $\M''$. We have the isomorphism $K''=K$ using the inductive definition of intersection complexes iterating open direct images and truncations, see \cite{BBD}. (Here we apply the last step of the inductive construction.) This implies that $\M''$ is a mixed Hodge module (that is, $H^k\M''\eq 0$ ($k\ne0$)), and its injective image in the mixed Hodge module $H^0(j_x)_*j_x^*\M$ is identified with the injective image of $\M$ in it, since this holds for the underlying $\Q$-complexes. (Note that $H^{\ssb}$ is the standard cohomology functor of the bounded derived category $D^b{\rm MHM(X)}$.) Thus $\M''$ in (2.1.7) can be replaced by $\M$.
\sk
The assertion (2.1.1) then follows from the standard estimates of weights for the pullback functor, see \cite[(4.5.2)]{mhm}. (Here it is also possible to use the ``classical" $t$-structure ${}^c\tau_{\les p}$ on the bounded derived category of mixed Hodge modules, see \cite[Remark 4.6,2]{mhm}.)
\msn
{\bf Remark~2.1b.} It does not seem necessarily easy to follow some arguments in \cite{ES}. For instance, the authors hire the theory of mixed Hodge modules {\it partially\1} in some places, although it does not seem quite clear whether the quoted assertions can really adapt to the situation they are considering, since they are performing a too complicated calculation of nearby cycles extending an old {\it double complex construction\1} in terms of logarithmic complexes and $\frac{\ddd f}{f}\wedge$ {\it without using filtered $\D$-modules\1} (see also \cite{ELY}).
Note also that the Hodge filtration can never be captured as in \cite[6.1.1]{ES} using a filtration in the abelian full subcategory of $D^b_c(X,\C_X)$ constructed in \cite{BBD}.
\msn
{\bf Remark~2.1c.} It seems that they have recently written another paper which is similar to the above one. It seems to claim that the decomposition theorem for the direct image of an intersection complex with coefficients in a polarizable variation of Hodge structure can be proved rather easily by employing Kashiwara's results on infinitesimal mixed Hodge structures {\it without\1} hiring $\D$-modules. It contains however certain {\it fundamental and conceptual defects\1} as follows.
\sk
First of all, the theory of ``perverse truncation of Hodge filtration" {\it never\1} exists. This can be seen in the case $\pi\,{:}\,V'\tos V$ is a {\it smooth projective\1} morphism onto a {\it smooth\1} variety in the proof of a theorem in the paper, where the Whitney stratification of $V$ has only {\it one stratum,} so the {\it perverse\1} truncation coincides with the {\it usual\1} one defined by using the {\it kernels\1} of differential up to a shift. ({\it Never say that this case is not needed for the proof. If the argument works in the birational case, it should do also in that case.}) Here one has to {\it separate\1} the differential of the {\it relative de~Rham complex\1} $\Om_{V'/V}^{\ssb}$ from that of the {\it de~Rham complex on the base space\1} $V$ via the equality
$$\RR\pi_*\C_X[d_{V'}]\eq{\rm DR}_V(\RR\pi_*\Om_{V'/V}^{\ssb}[d_{V'}{-}d_V]),$$
assuming $V'\eq W{\times}V$ for simplicity. Indeed, in the case $\pi\eq id$ for instance, we see that
$$\Hc^k(F^p\Om_V^{\ssb}[d_V])\eq\begin{cases}{\rm Ker}\,\ddd^p\,({\subset}\,\Om_V^p)&\h{if}\,\,\,k\eq p\mi d_V,\\ \,0&\h{if}\,\,\,k\,{\neq}\,p\mi d_V,\end{cases}$$
where $F^p\Om_V^{\ssb}\eq\sigma_{\ges p}\Om_V^{\ssb}\eq\Om_V^{\ssb\ges p}$. This is far from what we expected. The above ``separation" {\it cannot\1} be generalized without employing the complex of {\it induced\1} $D$-modules associated with a {\it filtered differential complex.}
\sk
This is the most crucial error in the paper. This ``filtered truncation" is used everywhere in the paper. It seems impossible to fix this without hiring ``filtered $D$-modules".
\sk
The second major problem seems to be that the ``nearby and vanishing cycle functors" are defined {\it only\1} in the ``normal crossing" case by employing ``logarithmic differential forms", since the theory of such forms is {\it not\1} sufficiently developed in the general case. It is then impossible to prove the ``commutativity" of these functors with the direct image functor ``in a compatible way with the Hodge filtration" in general. However, this commutativity is ``implicitly" utilized {\it everywhere in the paper.} For the classical case of a {\it one-parameter degeneration\1} studied by Steenbrink and others, one can apply the ``semicontinuity argument" instead of the commutativity in order to prove the ``freeness" of the Hodge filtration on the Deligne extension. However there is {\it no\1} generalization of such a simple argument, and it is quite hard to get this commutativity without showing the ``stability" of the $V$\!-filtration and the ``bistrictness" of $F,V$ on the direct image complex.
\sk
In the normal crossing case, the ``combinatorial description of filtered regular holonomic $\D$-modules of normal crossing type" is employed under the name of ``de~Rham family". For instance, the ``logarithmic intersection complex" is a subcomplex of the logarithmic complex, and these are subcomplexes of the ``de~Rham complexes" of the intersection and meromorphic $\D$-modules associated with a variation of Hodge structure. These subcomplexes are obtained by using the ``truncation" by the ``$V$\!-filtrations $V_{(i)}^{\ssb}$ along local irreducible components". This truncation induces filtered quasi-isomorphisms. This procedure can be extended to the case of the ``localization along a subdivisor" of the logarithmic intersection complex. We get the ``combinatorial description" by restricting to the highest degree forms and dividing it by $\msum_i\1V^{>1}_{(i)}$, since it is already truncated by $\mcap_i\1V^0_{(i)}$. The Hodge filtration $F$ is obtained by {\it iterating\1} the $\Gr_{V_{(i)}}^{\ssb}$, where the ``compatibility" of the $(n{+}1)$-filtrations $F,V_{(1)},\dots,V_{(n)}$ is indispensable for the well-definedness.
\sk
Here it should be noted that this ``combinatorial description" {\it never\1} implies an ``equivalence of categories". Indeed, the associated filtered $\C$-vector spaces give only the {\it approximating\1} ``nilpotent orbits", and a lot of information is lost about the ``Hodge filtration" by passing to the limit. In particular, a decomposition in the category of infinitesimal mixed Hodge structures {\it never\1} guarantees the ``compatibility" of the decomposition of a filtered sheaf complex with the ``Hodge filtration", and some more argument is required. It is however quite difficult to prove it without using the ``compatibility" of the $(n{+}1)$-filtrations $F,V_{(1)},\dots,V_{(n)}$.
\sk
Note also that ``Kashiwara's combinatorial description" explained in the paper is due to Galligo, Granger, and Maisonobe (although it may be essentially his work; it is rather interesting that their paper did not appear from ASENS as is written in the references of his paper). This is however a ``topological" description using the local cohomologies along the ``closed half line" $\R_{\ges0}\sst\C$ in the one-dimensional case, and is essentially the same as applying Verdier's extension theorem for monodromical sheaves inductively, where the compositions of nearby or vanishing cycle functors along local coordinates are used. The morphisms between them are induced by the canonical and (topological) variation morphisms. However these are {\it never\1} good for Hodge theory, since their compositions are the ``topological variations" $T_i{-}\1{\rm id}$ instead of $N_i\eq\log T_i$ in the unipotent monodromy case. So they {\it never\1} give ``morphisms of mixed Hodge structures". One should absolutely use the ``$D$-module-theoretic description" explained above. In the non-quasi-unipotent monodromy case, one has to choose however a ``section" of the exponential map $e^{2\pi i*}\,{:}\,\C\onto\C^*$ for this.
\sk
Kashiwara's theorem in the paper asserts {\it only\1} that a mixed Hodge structure endowed with commuting nilpotent endomorphisms and a polarization defined by using a residue is a ``nilpotent orbit". One has to prove that it coincides with the ``combinatorial description" of the ``nearby cycle complex" of the intermediate direct image of a polarizable variation of Hodge structure, which should be isomorphic to $\it\Psi^*\!I\!\Lc$, in order to prove the decomposition theorem, and some calculations about the nearby cycle and dual functors are needed for this. However, these do not seem to be done in the paper. (For instance, the Hodge filtration never appears in SGA7 XIV.) A standard argument for the Hodge filtration uses the cokernel of the inclusion $j_!\M[N]\,{\into}\,j_*\M[N]$ in the category of inductive limits of regular holonomic $D$-modules of normal crossing type, where $\M$ is the localization along $f^{-1}(0)$ of the intersection $\D$-module and $j_!\eq\DD\ssc j_*\ssc\DD$.
\sk
One problem is that the definition of the nearby cycle complex $\it\Psi^*\!I\!\Lc$ in the paper is not easy to understand. Indeed, the information of the finite-dimensional complex $\it\Psi^*\!IL$ associated with the ``de~Rham family" in the paper is not enough for this, and we need the ``de~Rham family" itself, that is, ``the combinatorial description of the nearby cycles" as a filtered regular holonomic $D$-module of normal crossing type (which is known in $D$-module theory). Note that the ``tilde embedding" depends on local coordinates, hence the well-definedness is not quite trivial. Also the reason for the definition $A_i\defs m_iN{-}N_i$ is not well explained. This is closely related to the double complex construction using $\ddd f\!/\!f\wedge$. It comes also from the relation $x_i\dd_{x_i}\delta(t{-}f)\eq{-}m_i\dd_tt\delta(t{-}f)$ in $D$-module theory.
Concerning the polarization defined by using the residue, one has to show that it coincides with the polarization {\it induced by the nearby cycle functor.} This is quite nontrivial.
These arguments are highly sophisticated, and it does not seem very clear whether they are really written by the authors.
\sk
In the paper the finite-dimensional filtered complex $\it\Psi^*\!IL$ is defined to be quasi-isomorphic to the ``stalk" at $x$ of the nearby cycle complex $\!\it\Psi^*\!I\!\Lc$ ``forgetting the Hodge filtration $F$". Note first that the ``stalkwise" decomposition of a graded piece of $\it\Psi^*\!I\!\Lc$ at each point {\it never\1} gives a decomposition of the graded piece into a direct sum of intersection complexes unless there is a certain ``relation" between the decompositions at various points, since a lot of information is lost by taking the restriction functor $i_x^*$ with $i_x\,{:}\,\{x\}\,{\into}\,X$ the natural inclusion. So the ``de~Rham family" itself, that is, the ``combinatorial description" as a ``filtered regular holonomic $D$-module of normal crossing type" is needed.
\sk
Second the ``relation" between the Hodge filtrations $F$ on $\it\Psi^*\!IL$ and $\it\Psi^*\!I\!\Lc$ is quite unclear. The former {\it never\1} coincides with the ``stalk" of the Hodge filtration on the nearby cycle complex $\!\it\Psi^*\!I\!\Lc$ at $x$, since the latter is {\it infinite-dimensional.} It should coincide with the Hodge filtration obtained by iterating the graded quotients of $V_{(i)}$. Then one has to show the ``compatibility" of the $(n{+}2)$-filtrations $F,W(N),V_{(1)},\dots,V_{(n)}$ in order to prove that the decomposition of the $\Gr^{W(N)}_{\ssb}\it\Psi^*\!I\!\Lc$ is ``compatible with the Hodge filtration".
\sk
Concerning ``the Verdier extension theorem with Hodge filtration", this is {\it never\1} formulated in the paper. It seems quite difficult to state and prove it without using the $V$\!-filtration due to the {\it lack\1} of a theory of {\it logarithmic\1} forms in the {\it non-normal crossing\1} case. (There is a so-called ``logarithmic comparison theorem" which is proved only for certain cases where the variation of Hodge structure must be locally constant.)
\sk
As a conclusion, the proofs of  the ``compatibilities with Hodge filtration" are {\it neglected\1} ``systematically". It seems quite difficult to fix the problems without using $D$-modules and $V$\!-filtrations, although the resulting theory could become quite similar to that of mixed Hodge modules.
\msn
{\bf Remark~2.1d.} It is quite surprising that the authors of the paper are still continuing their {\it fantasy\1} using their ``magic $t$-structure" although it is absolutely impossible to control the Hodge filtration by employing only the $t$-structure.
\sk
Let $K$ be a bounded constructible complex on a variety $V$ endowed with two filtrations $F_1,F_2$ in the usual sense (not necessarily by constructible complexes). They {\it correctly\1} claim that there is a bifiltered bounded complex $(K';F'_1,F'_2)$ together with a bifiltered quasi-isomorphism $(K;F_1,F_2)\to(K';F'_1,F'_2)$ such that the truncations $^p\tau_{\le j}K'$ ($j\ins\Z)$ can be {\it represented\1} by subcomplexes of $K'$. This assertion however never implies the {\it well-definedness of the induced filtrations\1} $F'_1,F'_2$ on $^p\Hc^jK'$. Indeed, these heavily depend on the {\it representative\1} of the truncation $^p\tau_{\le j}K$ ($j\ins\Z)$.
\sk
This phenomenon can be observed typically in the following case: Assume $V$ is smooth and $K$ is the usual de Rham complex $\Om_V^{\ssb}$ with $F_1$ the Hodge filtration $F$ defined by $F^k\defs\sigma^{\ge k}$ and $F_2$ the weight filtration $W$ defined by $\Gr^W_i\!K\eq0$ ($i\nes 0$). We have $^p\Hc^j K\eq0$ for $j\ne d_V\defs\dim V$ as well known. We can nevertheless apply the truncation $^p\tau_{\le d_V}$ to this complex. The proof of the above assertion uses an addition of the mapping cone of an identity morphism, but this can be skipped when the representative of $^p\tau_{\le d_V}K$ is given by a subcomplex (more precisely, the result does not change; we cannot avoid taking the intersection anyway in order that the morphism preserve the filtrations). This applies to the case $^p\tau_{\le d_V}K$ is {\it represented by\1} $\tau_{\le 0}\Om^{\ssb}_V\eq\C_V$ following the definition of $^p\tau$ and using Poincar\'e lemma at least on a non-empty Zariski-open subset of $V$. This representative of $^p\tau_{\le d_V}K$ is however very bad for the Hodge filtration $F^k\defs\sigma^{\ge k}$ if $k\sgt 0$, since their intersection vanishes. We may take however $K$ itself as the representative of $^p\tau_{\le d_V}K$. So in the general case {\it very good representative of the truncation\1} $^p\tau_{\le j}$ must be chosen. This is however an {\it extremely difficult problem\1} (without using $\D$-modules).
\sk
The authors claim that the decomposition in the decomposition theorem for the direct image of the intersection complex by a projective morphism (for instance a desingularization) is {\it compatible with the Hodge filtration.} This should however imply the {\it strictness\1} of the direct image of the corresponding filtered $\D$-module, and is never trivial. It should impose an additional difficulty to treat the Hodge filtration on intersection complexes in the paper. Indeed, the {\it strictness\1} should hold for the direct factor of the direct image supported on points. So its proof must be given (including the case $f(X)$ is a point).
The compatibility of the decomposition with the Hodge filtration should immediately imply for instance the pure Hodge structure on the intersection cohomology of a projective variety with constant coefficients by using only the classical Hodge theory together with the decomposition theorem without Hodge filtration. The latter assertion is never trivial even in the isolated singularity case as seen in Navarro's paper.
\sk
One of the main problems of the paper is that it is quite unclear in which category the decomposition isomorphism is proved. {\it This must be clarified absolutely.} One may apply Deligne's argument to spectral objects in the sense of Verdier, but never to ``filtered" ones in the sense of the paper. They should  also have to show at least the ``strictness" of the Hodge filtration on the direct factors {\it supported on points.} It is completely unclear how these problems are solved by employing their ``magic $t$-structure".
\sk
As for the proof of the decomposition theorem for the direct image of the intersection complex, it has nothing to do with Gabber in the {\it constant\1} coefficient case. His contribution to this theory is the proof of the {\it local purity\1} of intersection complexes with non-constant coefficients. (The stability of pure complexes by the direct image under a proper morphism was shown in [Weil II], and the decomposition of pure complexes {\it forgetting the Frobenius action\1} is proved in [BBD].)
\sk
Concerning the criterion for the decomposition using the decomposition of the vanishing cycles $\varphi_{f,1}$ as the direct sum of ${\rm Im\,can}$ and ${\rm Ker\,Var}$, this is {\it never\1} written in Verdier's papers even though it is a direct corollary of Verdier's extension theorem {\it as long as the Hodge filtration is forgotten,} provided that one notices it. When it was stated at the first time without mentioning Verdier, nobody (including him) complained. Note that this criterion implies for instance that the decomposition holds with $\Q$-coefficients in the projective case applying Deligne's argument to spectral objects. This should have been mentioned in his paper if he knew it, since it is a very good application of the theorem. (Recall that the decomposition was proved only with $\C$-coefficients in [BBD].)
\sk
In order to apply Lemma 5.2.15 in [MHP], one must be very careful about the sign. Indeed, if the sign is wrong, the roles of $H$ and $H'$ are exchanged. It is however quite unclear how this problem is solved in the paper. It is known that the duality isomorphism contains some sign problem, so it is better to use the direct image of the corresponding self-pairing as in [MHP], see also 2103.04836. But this kind of self-pairing does not seem to be constructed in the paper.
\sk
It seems very difficult to prove the compatibility of the decomposition with the Hodge filtration $F$ {\it without\1} using the $V$-filtration indexed by $\Q$. Indeed, it might occur that the filtration $F$ is {\it not\1} compatible with the decomposition although it becomes compatible by {\it taking the unipotent vanishing cycle functor\1} $\varphi_{f,1}$. Here the definition of {\rm Var} using the action of $t$ together with the $V$-filtration is quite useful in order to prove this compatibility in the case it holds {\it after taking $\varphi_{f,1}$ or equivalently\1} $\Gr_V^0$. (For instance we have the surjection ${\rm Ker}\,t\cap V^0\onto{\rm Ker}\,\Gr_V^0t$, which is very important for the study of the filtration $F$.)
\sk
As for the ``topological" $N'\defs T_u\mi 1\eq\msum_{k\ge 1}N^k\!/k!$, note that this never preserves the Hodge filtration {\it up to the shift by\1} 1, hence it never gives a morphism of mixed Hodge structure of type $(-1,-1)$ in general.
\sk
Concerning the definition of admissible ``PVMHS" (which is {\it not\1} polarizable VMHS) in the second part, the {\it existence\1} of the Hodge filtration is assumed {\it only on $\psi_f$ without using the Deligne extension or $V$-filtration.} This formulation is however {\it totally insufficient\1} in the classical case of a variation of mixed Hodge structure on a curve studied forty years ago by one of the authors. Here one has to define the limit Hodge filtration on $\psi_f$ using the Deligne extensions and assume that the induced filtration on each graded piece of the weight filtration $W$ coincide with the limit Hodge filtration for the graded piece of $W$. This is a typical example of lots of fake arguments in the paper.
\sk
In the case we have a polarizable variation of Hodge structure on a smooth open dense subvariety of $X$ together with a projective morphism $f\,{:}\,X\tos Y$ and a projective birational morphism $\pi\,{:}\,X'\tos X$, it seems very difficult to deduce the hard Lefschetz property for the direct image by $f$ of the the intermediate direct image of the polarizable variation from that for the direct image by $f\ssc\pi$. There is a counterexample in the Hodge setting to some ``argument" in the twistor theory, see 2204.09026. This ``argument" is quite long and complicated (involving the construction of Chern classes), and is never a ``linear algebra" as is noted in 2401.09544. It is very difficult to find out the exact place where an error occurred.
\sk
Verdier's extension theorem is called ``Verdier's classification" in the paper (perhaps in order to give some illusion). However it is an {\it equivalence of categories,} and never provides a classification. Its formulation and proof in the paper are quite different from those in Verdier's paper, where the reduction to the {\it monodromical\1} case is used employing Verdier specialization together with an argument due to Deligne. Moreover the nearby and vanishing cycle functors were {\it not\1} restricted to their {\it unipotent\1} monodromy part. Indeed, the restriction to the unipotent monodromy part of the nearby and vanishing cycle functors was never written in his paper. This is only after the theory of mixed Hodge modules where the $V$-filtration {\it indexed by\1} $\Q$ is essential.
\sk
Verdier's extension ``Lemma" is proved in the paper {\it without\1} employing the specialization although the authors say that Verdier's proof is {\it adapted\1} to their situation. In this case they should explain why the specialization had to be used in his paper, and how the problem is avoided in their paper. It is quite interesting that the answers seem to be available by looking for instance at the proof of the ``essential surjectivity" in Verdier's extension ``Lemma" in the paper. Here a complex $F_2$ is defined to be the mapping cone of some extension class. However there is no explanation about the construction of this class, which seems very difficult to specify even by using the Leray type spectral sequence for higher extension groups. They demonstrate only trivial properties which can hold even for the mapping cone of the {\it trivial\1} (that is, splitting) extension class, and the proof suddenly ends without giving any explanation about the last crucial isomorphism, where $F$ seems to be $F_2$ (and the proof contains some serious typos). It does not seem that the theorem can be proved in this way, and this seems to give the answers to the above questions. This is another typical example of lots of fake arguments in the paper.
\sk
As for the compatibility of Verdier's extension theorem with the Hodge filtration, it is still left completely untouched. The assertion does not seem to be proved correctly without using the $V$-filtration of Kashiwara and Malgrange {\it indexed by\1} $\Q$ as is stressed repeatedly.
The Hodge filtration on the nearby or vanishing cycle functor in the {\it non-normal crossing case\1} is always defined in the paper by using the direct image by a proper morphism of the nearby or vanishing cycle functor {\it in the normal crossing case.} So the arguments suppose from the beginning an ``imaginary commutativity" of the nearby and vanishing cycle functors with direct images. This may be called really a ``house built on sand". Indeed, for the geometric proof of Schmid theorem by one of the authors almost fifty years ago, the {\it freeness\1} of the higher direct image sheaves of the relative logarithmic differential forms was crucial, however this part is {\it completely forgotten\1} and the argument is {\it simplified in a revolutionary way\1} in the paper, see also the remark on the definition of admissible ``PVMHS".
\sk
As a conclusion, the authors pay no serious attention to the problem of the compatibility with the Hodge filtration supplying lots of fake arguments, and never try to understand properly the difficulty of the problems related to the Hodge filtration. Their theory is still a ``fantasy" or a ``house built on sand" just as before

\ms
{\smaller\smaller RIMS Kyoto University, Kyoto 606-8502 Japan}
\end{document}